\documentclass[12pt]{article}

\usepackage{amssymb}
\usepackage{amsmath}
\usepackage{graphicx}
\usepackage{hyperref}

\usepackage{textgreek}

\numberwithin{equation}{section}

\begin{document}

\newcommand{\bd}{\begin{displaymath}}
\newcommand{\ed}{\end{displaymath}}
\newcommand{\ds}{\displaystyle}
\newcommand{\bp}{\underline{\bf Proof}:\ }
\newcommand{\ep}{{\hfill $\Box$}\\ }
\newcommand{\be}{\begin{equation}}
\newcommand{\ee}{\end{equation}}
\newcommand{\ba}{\begin{array}}
\newcommand{\ea}{\end{array}}
\newcommand{\bea}{\begin{eqnarray}}
\newcommand{\eea}{\end{eqnarray}}
\newcommand{\nt}{\noindent}

\newtheorem{0}{DEFINITION}[section]
\newtheorem{1}{LEMMA}[section]
\newtheorem{2}{THEOREM}[section]
\newtheorem{3}{COROLLARY}[section]
\newtheorem{4}{PROPOSITION}[section]
\newtheorem{5}{REMARK}[section]
\newtheorem{6}{EXAMPLE}[section]
\newtheorem{7}{ALGORITHM}[section]
\newtheorem{8}{CONJECTURE}[section]

\title{Mean First Passage Times and Fundamental Tensors of Higher Order Markov Chains}
\author{
Jianhong Xu\thanks{School of Mathematical and Statistical Sciences, Southern Illinois University Carbondale, Carbondale, IL 62901, USA. Email: \texttt{jhxu@siu.edu}} 
% \and  \thanks{}
% \and  \thanks{}
}

\maketitle

\begin{abstract}  
The mean first passage times are among the most critical characteristics of a Markov chain. In this paper, we focus on the scenario in which one or more states of a higher order ergodic Markov chain are modified to be absorbing. We prove that the resulting chain has to be absorbing. For a higher order absorbing Markov chain, we prove that the equation its fundamental tensor satisfies must be nonsingular and provide a MATLAB function {\tt fund} for solving the equation. Besides, we connect each horizontal slice of the mean first passage time tensor with a fundamental tensor obtained when one state of a higher order ergodic Markov chain is modified to be absorbing, which also leads to a tensor series representation for selected mean first passage times.
\end{abstract}

\nt {\bf Keywords}: higher order Markov chain, transition tensor, tensor multiplication, ergodicity, mean first passage times, higher order absorbing Markov chain, canonical form, fundamental tensor

\nt {\bf AMS Subject Classification}: 15A69, 15B51, 60A05, 60J10, 60J99 

\section{Introduction}
\label{intro}
\setcounter{equation}{0}

Markov chains are one of the most important classes of stochastic processes to model systems whose states evolve over time according to certain probabilistic transition rules. Classical first order Markov chains suppose that the future state depends solely on the current state, not on any preceding states. For many real-world systems, however, the future state depends on not only the current state but also a couple of preceding states. Higher order Markov chains, an extension to the first order case, serve as a natural framework for modeling such phenomena.

Specifically, the notion of $(m-1)$th order Markov chains is as follows. Let $m, n \ge 2$ be integers. Consider a stochastic process
$X=\{X_t : t=1, 2, \ldots\}$. For each $t$, the random variable $X_t$ takes values in $S=\{1, 2, \ldots, n\}$. Then, the process $X$ is called an $(m-1)$th order Markov chain if for all $t \ge m-1$ and $i_1, i_2, \ldots, i_m \in S$, we have 
\be
\label{markov}
\begin{array}{l}
\Pr(X_{t+1}=i_1 | X_t=i_2,\ldots,X_{t-m+2}=i_m,\ldots,X_1=i_{t+1})\\
=\Pr(X_{t+1}=i_1 | X_t=i_2,\ldots,X_{t-m+2}=i_m).
\end{array}
\ee
The set $S$ is called the state space. When $t$ is interpreted as time and $X_t$ as the state of a system, (\ref{markov}) shows precisely that the future state $X_{t+1}$ hinges upon the current state $X_t$ and $m-2$ preceding states $X_{t-1}$ through $X_{t-m+2}$. In particular, the special case when $m=2$ gives rise to a classical first order Markov chain and the case when $m \ge 3$ yields a higher order Markov chain.

In addition, a Markov chain is called homogeneous if the probabilities in (\ref{markov}) are independent of $t$, i.e.,
$$\Pr(X_{t+1}=i_1 | X_t=i_2,\ldots,X_{t-m+2}=i_m)=\Pr(X_m=i_1 | X_{m-1}=i_2,\ldots,X_1=i_m)$$
for any $t \ge m$. Accordingly, each probability there can be denoted as
$$\Pr(X_{t+1}=i_1 | X_t=i_2,\ldots,X_{t-m+2}=i_m)=p_{i_1i_2\ldots i_m}, ~t \ge m-1,$$
and called the transition probability from states $(i_2, \ldots, i_m)$ to $i_1$. Its subscript $i_1i_2\ldots i_m$ is sometimes referred to as a multi-index of length $m$. The $m$th order, $n$ dimensional array ${\cal P}=[p_{i_1i_2\ldots i_m}]$ is called the transition tensor. Since $0 \le p_{i_1i_2\ldots i_m} \le 1$ for any $i_1, i_2, \ldots, i_m \in S$ and $\ds \sum_{i_1 \in S}p_{i_1i_2\ldots i_m}=1$ for any $i_2, \ldots, i_m \in S$, $\cal P$ is said to be a stochastic tensor.

Throughout this work, we are concerned only with homogeneous Markov chains. For brevity, we shall refer to them simply as chains from now on.

For first order chains, a rich, well-established theory has long been available; see the classic treatises \cite{Doo,Hun,Ios,KS}. With recent advances in tensor theory and computations, see, for example, \cite{CW, DW, KB, MSL, QCC, QL}, along with many emerging applications \cite{Bae, Bur, Fle, Ho, Isl, Kwa, Lan, LLZ, M, San, Xio, YJK}, higher order chains have drawn much attention within the research community over the past decade or so \cite{CZ, CPZ, G, GLY, HQ, HWX, HX24a, HX24b, HX26a, HX26b, LZ, LN, WC, Xu26a, Xu26b}.

In particular, the mean first passage times of a higher order chain have been studied in \cite{HWX, HX24a}. Given the $(m-1)$th order chain $X$ on state space $S=\{1, 2, \ldots, n\}$ as specified earlier, the first passage time random variable $\eta_{i_1i_2\ldots i_m}$ from states $(i_2, \ldots, i_m)$ to $i_1$ is defined by 
\be
\label{eta}
\eta_{i_1i_2\ldots i_m}=\min\{k \ge 1 : X_{m+k-1}=i_1 | X_{m-1}=i_2, \ldots, X_1=i_m\}.
\ee
The mean first passage time from states $(i_2, \ldots, i_m)$ to $i_1$ follows as 
$$\mu_{i_1i_2\ldots i_m}={\rm E}(\eta_{i_1i_2\ldots i_m})=\sum_{k=1}^\infty k\Pr(\eta_{i_1i_2\ldots i_m}=k).$$
The $m$th order, $n$ dimensional tensor $\mu=[\mu_{i_1i_2\ldots i_m}]$ is called the mean first passage time tensor.

Such mean first passage times are crucial because they can characterize the short-term behavior of a chain. According to S. Redner \cite{Red}, ``first passge underlies many stochastic processes in which the event, such as ... a chemical reaction, the firing of a neuron, or the triggering of a stock option, relies on a variable reaching a specified value for the first time.''

In the meantime, the notion of higher order absorbing chains has recently been proposed in \cite{HX26b}. Detail in this regard will be given in the next section. For the first order case, the relationship between the mean first passage times and the expected times to absorption is known \cite{KS}. It is natural, therefore, to ask whether such relationship can be generalized to the higher order case.

It should be pointed out that as we go from the first order to the higher order cases, quite a few differences arise \cite{HX24a, HX24b, HX26a, Xu26a}. Although it is well known that any higher order chain can be associated with a first order chain \cite{Doo,Hun}, not all the problems regarding a higher order chain can be addressed by turning to its associated first order chain \cite{Xu26a}. While it comes to higher order absorbing chains, we are facing the same situation \cite{HX26b}.

In this work, we shall investigate the relationship between the mean first passage time tensor of an $(m-1)$th order chain and the fundamental tensor of a relevant $(m-1)$th order absorbing chain. As in \cite{HWX, HX24a, HX26b}, our focus will be on the higher order case, i.e., $m \ge 3$. Our results, of course, apply to the first order case as well by allowing $m=2$.

Specifically, our main goals can be summarized as:
\begin{itemize}
\item{We shall show that if one or more states of an $(m-1)$th order ergodic chain $X$ are made absorbing, then the resulting chain $\tilde X$ is absorbing.}
\item{We shall develop several useful results regarding the fundamental tensor of the chain $\tilde X$, including the case of a general higher order absorbing chain, especially regarding the nonsingularity of a tensor equation governing the fundamental tensor.}
\item{We shall give a useful MATLAB function {\tt fund} for computing the fundamental tensor of an absorbing chain of any order.}
\item{We shall establish the connection between the mean first passage times of the chain $X$ to a state $i$ and the fundamental tensor of the chain $\tilde X$ with state $i$ being a single absorbing state, which also leads to a tensor series representation for selected mean first passage times.}
\item{We shall also provide necessary examples to illustrate our results.}
\end{itemize}

Our presentation will be organized as follows. In Section \ref{back}, some necessary background material is cited. Our main results are developed in Section \ref{main}, which are accompanied by several illustrative examples. Some concluding remarks are given at the end in Section \ref{concl}.

Throughout this paper, $\mathbb R^{n_1 \times n_2 \times \cdots \times n_m}$ will be used to denote the space of $m$th order tensors on the real field $\mathbb R$ of size $n_1 \times n_2 \times \cdots \times n_m$, with $n_i \ge 1$ being integers for $1 \le i \le m$. 

\section{Background Material}
\label{back}
\setcounter{equation}{0}

Let us start with a tensor $\boxtimes$ product and its resulting tensor power in \cite{HX26b}.

\begin{0}
\label{bprod}
Given $m$th order tensors ${\cal A}=[a_{i_1i_2\ldots i_m}] \in \mathbb R^{u \times w \times n \times \cdots \times n}$ and ${\cal B}=[b_{i_1i_2\ldots i_m}] \in \mathbb R^{w \times v \times n \times \cdots \times n}$, where $v \le n$, ${\cal C}={\cal A}\boxtimes {\cal B}=[c_{i_1i_2\ldots i_m}] \in \mathbb R^{u \times v \times n \times \cdots \times n}$ is defined by 
$$c_{i_1i_2\ldots i_m}=\sum_{j=1}^w a_{i_1ji_2\ldots i_{m-1}}b_{ji_2\ldots i_m}$$
for any $1 \le i_1 \le u$, $1 \le i_2 \le v$, and $1 \le i_3, \ldots, i_m \le n$.
\end{0}

The reason for restricting $v \le n$ can be seen from the third index of the factor $a_{i_1ji_2\ldots i_{m-1}}$. Such a tensor $\boxtimes$ product is distributive but, in general, it is not associate when $m \ge 3$. Clearly, when $m=2$, it reduces to the regular matrix multiplication. The above definition leads to:  

\begin{0}
\label{bpow}
Given $m$th order tensor ${\cal A}=[a_{i_1i_2\ldots i_m}] \in \mathbb R^{u \times u \times n \times \cdots \times n}$, where $u \le n$, the $k$th power of $\cal A$ is defined recursively by 
$${\cal A}^1={\cal A}, ~{\cal A}^{k+1}={\cal A}^k \boxtimes {\cal A}, ~k=1, 2, \ldots.$$
By convention, ${\cal A}^0$ is the identity tensor ${\cal I}=[\delta_{i_1i_2\ldots i_m}] \in \mathbb R^{u \times u \times n \times \cdots \times n}$, where for any $1 \le i_1, i_2 \le u$ and $1\le i_3, \ldots, i_m \le n$,
$$\delta_{i_1i_2\ldots i_m}=\delta_{i_1i_2}=\left\{\begin{array}{cl}
1, & i_1=i_2;\\
0, & {\rm otherwise}.
\end{array}\right.$$
\end{0}

It is easy to verify that $\cal I$ is a left identity tensor, i.e., ${\cal I} \boxtimes {\cal B}={\cal B}$ for any $m$th order tensor ${\cal B} \in \mathbb R^{u \times v \times n \times \cdots \times n}$ with $v \le n$; but, in general, ${\cal A}\boxtimes {\cal I} \ne {\cal A}$ for an $m$th order tensor ${\cal A} \in \mathbb R^{v \times u \times n \times \cdots \times n}$ with $u \le n$.

The preceding definitions extend the product, power, and identity tensor in \cite{HX24a, HX26a} for the special case $u=w=v=n$. All these, including the case in \cite{HX24a, HX26a}, have been implemented in \cite{Xu26c} using MATLAB as {\tt bprod}, {\tt bpow}, and {\tt eyet}, respectively.

Continuing, the notion of a diagonal tensor can be stated as:
\begin{0}
For an $m$th order tensor ${\cal A}=[a_{i_1i_2\ldots i_m}] \in \mathbb R^{u \times u \times n \times \cdots \times n}$, the diagonal tensor obtained from $\cal A$, denoted by ${\cal A}_d=[a^{[d]}_{i_1i_2\ldots i_m}]$, is defined as: for any $1 \le i_1, i_2 \le u$ and $1 \le i_3, \ldots, i_m \le n$, 
$$a^{[d]}_{i_1i_2i_3\ldots i_m}=a_{i_1i_2i_3\ldots i_m}\delta_{i_1i_2}=\left\{\begin{array}{cl}
a_{i_1i_1i_3\ldots i_m}, & i_1=i_2;\\
0, & {\rm otherwise}.
\end{array}\right.$$
\end{0}

Returning to the $(m-1)$th order chain $X$ with $m$th order, $n$ dimensional transition tensor ${\cal P}=[p_{i_1i_2\ldots i_m}]$ and state space $S=\{1, 2, \ldots, n\}$ as laid out in the introductory section, we first quote the technical lemma below. 

\begin{1}[\cite{HX24a}]
\label{kstep}
Let ${\cal P}^k=[p^{(k)}_{i_1i_2\ldots i_m}]$, $k=1, 2, \ldots$. Then, $p^{(k)}_{i_1i_2\ldots i_m}$ is the $k$-step transition probability from states $(i_2, \ldots, i_m)$ to $i_1$, i.e.,
$$p^{(k)}_{i_1i_2\ldots i_m}=\Pr(X_{m+k-1}=i_1 | X_{m-1}=i_2, \ldots, X_1=i_m),$$
for any $i_1, i_2, \ldots, i_m \in S$.
\end{1}

Consequently, ${\cal P}^k$ is called the $k$-step transition tensor. One usage of ${\cal P}^k$ is to specify the ergodicity of the chain $X$.
\begin{0}[\cite{HWX, HX24a}]
\label{ergo}
The chain $X$ is said to be ergodic if for any $i_1, i_2, \ldots, i_m \in S$, there exists $k \ge 1$, which may depend on $i_1, i_2, \ldots, i_m$, such that  $p^{(k)}_{i_1i_2\ldots i_m}>0$.
\end{0}

Regarding the mean first passage time tensor $\mu$ of the chain $X$, we have the following fundamental result. 
\begin{1}[\cite{HX24a}]
\label{mfpt}
If the chain $X$ is ergodic, then its mean first passage times are all finite and satisfy 
\be
\label{mfpteqn}
\mu_{i_1i_2\ldots i_m}=1+\sum_{j \in S, j \ne i_1}\mu_{i_1ji_2\ldots i_{m-1}}p_{ji_2\ldots i_m}
\ee
for any $i_1, i_2, \ldots, i_m \in S$. In tensor form, (\ref{mfpteqn}) can be written as 
$$\mu={\cal E}+(\mu-\mu_d)\boxtimes {\cal P},$$
where $\cal E$ is the tensor of all ones and of the same size as $\cal P$, and where $\mu_d$ is the diagonal tensor obtained from $\mu$. Moreover, as a linear system, (\ref{mfpteqn}) is nonsingular if and only if the chain $X$ is ergodic.
\end{1}

Next, an absorbing state of the chain $X$ can be described as: 
\begin{0}[\cite{HX24b}]
A state $i$ of the chain $X$ is said to be absorbing if $p_{iii_3\ldots i_m}=1$ for any $i_3, \ldots, i_m \in S$. 
\end{0}

Here is a direct consequence of a state being absorbing: 
\begin{1}[\cite{HX24b}]
\label{absopr}
Let $i$ be an absorbing state of the chain $X$. Then, for any $k \ge 1$ and $i_3, \ldots, i_m \in S$, 
$p^{(k)}_{iii_3\ldots i_m}=1$, i.e., 
\be
\label{absoeqn}
p^{(k)}_{i_1ii_3\ldots i_m}=\delta_{i_1i}
\ee
for any $i_1, i_3, \ldots, i_m \in S$.
\end{1} 

By (\ref{absoeqn}), it is obvious that an ergodic chain does not have any absorbing state.

Finally, we recall the following concept of an $(m-1)$th order absorbing chain.

\begin{0}[\cite{HX26b}]
\label{abso}
The chain $X$ is called absorbing if it has at least one absorbing state and for every non-absorbing state $i_1$ and any $i_2, \ldots, i_m \in S$,
\be
\label{absolim}
\lim_{k \rightarrow \infty}p^{(k)}_{i_1i_2\ldots i_m}=0.
\ee
\end{0}

Observe that if we denote by $T$ the set of all non-absorbing states, then condition (\ref{absolim}) is equivalent to 
$$\lim_{k \rightarrow \infty} p^{(k)}_{i_1i_2i_3\ldots i_m}=0$$
for any $i_1, i_2 \in T$ and $i_3, \ldots, i_m \in S$, since (\ref{absoeqn}) implies $p^{(k)}_{i_1i_2i_3\ldots i_m}=0$ whenever $i_1 \in T$ and $i_2 \notin T$.

As shown in \cite{HX26b}, condition (\ref{absolim}) is crucial in guaranteeing that the non-absorbing states of the chain $X$ are all transient and the chain $X$ will eventually be absorbed in some absorbing state. Notice that, as in the first order case, the existence of absorbing states alone does not translate into a higher order chain being absorbing \cite{HX26b}.

In the remainder of this paper, an absorbing chain will be denoted as $\tilde X$ to distinguish it from the chain $X$, which will be assumed ergodic for dealing with the mean first passage times.

\section{Main Results}
\label{main}
\setcounter{equation}{0}

Recall that $X$ is an $(m-1)$th order chain with state space $S=\{1, 2, \ldots, n\}$ and $m$th order, $n$ dimensional transition tensor ${\cal P}=[p_{i_1i_2\ldots i_m}]$.

With the background material in the previous section, we are now ready to establish new results regarding the mean first passage times of the chain $X$ in relation to the fundamental tensor of a certain type of modified chain. Due to Lemma \ref{mfpt}, the chain $X$ will be assumed as ergodic from now on.

Our first question concerns a modified chain $\tilde X$ that is obtained from $X$ by making one or more states absorbing. Assume, without loss of generality, that these are the last $1 \le a < n$ states; otherwise, we can simply perform a simultaneous permutation on all the indices $i_1, i_2, \ldots, i_m$ to permute all the absorbing states to the last $a$ positions.

Specifically, let us partition the transition tensor $\cal P$ as 
\be
\label{part0}
{\cal P}=\left[\begin{array}{cc}
{\cal Q} & {\cal S}\\
{\cal R} & {\cal T}
\end{array}\right],
\ee
where $\cal Q$, $\cal R$, $\cal S$, and $\cal T$ are all $m$th order tensors such that ${\cal Q} \in \mathbb R^{t \times t \times n \times \cdots \times n}$, ${\cal R} \in \mathbb R^{a \times t \times n \times \cdots \times n}$, ${\cal S} \in \mathbb R^{t \times a \times n \times \cdots \times n}$, ${\cal T} \in \mathbb R^{a \times a \times n \times \cdots \times n}$, and $t=n-a$. By this partitioned form, we mean that each frontal slice of $\cal P$ can be expressed as 
$${\cal P}(:,:,i_3,\ldots,i_m)=\left[\begin{array}{cc}
{\cal Q}(:,:,i_3,\ldots,i_m) & {\cal S}(:,:,i_3,\ldots,i_m)\\
{\cal R}(:,:,i_3,\ldots,i_m) & {\cal T}(:,:,i_3,\ldots,i_m)
\end{array}\right].$$ 
For convenience, set $T=\{1, 2, \ldots, t\}$ and $A=\{t+1, t+2, \ldots, n\}$. Under this setting, the transition tensor of $\tilde X$ is given by 
\be
\label{part}
\tilde{\cal P}=[\tilde p_{i_1i_2\ldots i_m}]=\left[\begin{array}{cc}
{\cal Q} & {\cal O}\\
{\cal R} & {\cal I}
\end{array}\right],
\ee
where ${\cal O}$ is the tensor of all zeros and of the same size as $\cal S$, and where ${\cal I}$ is the identity tensor of the same size as $\cal T$.

Before formulating our main results, we need one more technical lemma.

\begin{1}[\cite{HX26b}]
Let $\tilde {\cal P}$ be given as in (\ref{part}). Then, for any $k=1, 2, \ldots$, 
\be
\label{partpkeqn}
\tilde{\cal P}^k=[\tilde{p}^{(k)}_{i_1i_2\ldots i_m}]=\left[\begin{array}{cc}
{\cal Q}^k & {\cal O}\\
{\cal R}_k & {\cal I}
\end{array}\right],
\ee
where ${\cal R}_k$ satisfies 
$${\cal R}_1={\cal R}, ~{\cal R}_{k+1}={\cal R}_k \boxtimes {\cal Q} + {\cal R}.$$
Equivalently, we have the following partitioned tensor $\boxtimes$ multiplication: 
$$\left[\begin{array}{cc}
{\cal Q}^k & {\cal O}\\
{\cal R}_k & {\cal I}
\end{array}\right] \boxtimes \left[\begin{array}{cc}
{\cal Q} & {\cal O}\\
{\cal R} & {\cal I}
\end{array}\right]=\left[\begin{array}{cc}
{\cal Q}^{k+1} & {\cal O}\\
{\cal R}_{k+1} & {\cal I}
\end{array}\right].$$ 
\end{1}

The above lemma is shown in \cite{HX26b} on a stochastic tensor $\tilde {\cal P}$ in the form (\ref{part}) without the assumption that the originating chain $X$ is ergodic. Actually, it remains valid for any $m$th order, $n$ dimensional tensor $\tilde {\cal P}$ with a zero block $\cal O$ and an identity block $\cal I$ being situated as in (\ref{part}). Notice, however, that if (\ref{part}) is permuted to 
$$\tilde {\cal P}=\left[\begin{array}{cc}
{\cal I} & {\cal R}\\
{\cal O} & {\cal Q}
\end{array}\right],$$
then a conclusion parallel to (\ref{partpkeqn}) fails to hold, i.e., 
$$\tilde{\cal P}^k \ne \left[\begin{array}{cc}
{\cal I} & {\cal R}_k\\
{\cal O} & {\cal Q}^k
\end{array}\right].$$
Incidentally, (\ref{part}) is called the canonical form in the case when the chain $\tilde X$ is absorbing \cite{HX26b}.

\begin{2}
\label{bound}
Given the $(m-1)$th order ergodic chain $X$, let $\tilde X$ be the modified chain obtained from $X$ by making $1 \le a < n$ states absorbing. We denote by $T$ the set of non-absorbing states and $\tilde {\cal P}=[\tilde p_{i_1i_2\ldots i_m}]$ the transition tensor of $\tilde X$. Then, there exist $0 \le \alpha < 1$ and $K \ge 1$ such that 
\be
\label{tpk}
\sum_{i_1 \in T}\tilde p^{(k)}_{i_1i_2\ldots i_m} \le \alpha^k
\ee
for all $i_2, \ldots, i_m \in S$ and $k \ge K$.
\end{2} 
\bp
Without loss of generality, we assume $\tilde {\cal P}$ as in the form (\ref{part}). Similar to (\ref{eta}), let $\eta_{Ai_2\ldots i_m}$ be the random variable of first passage time on $X$ from states $(i_2, \ldots, i_m)$ to $A=\{t+1, t+2, \ldots, n\}$, i.e.,
 $$\eta_{Ai_2\ldots i_m}=\min \{j \ge 1 : X_{m+j-1} \in A | X_{m-1}=i_2, \ldots, X_1=i_m\}.$$ 
Since $X$ is ergodic, by Definition \ref{ergo}, for any $i_1 \in A$ and $i_2, \ldots, i_m \in S$, there exists $k \ge 1$, which may depend on $i_1, i_2, \ldots, i_m$, such that $p^{(k)}_{i_1i_2\ldots i_m} > 0$. Let 
$$K=\max_{i_2, \ldots, i_m \in S}\min_{i_1 \in A} \{k \ge 1 : p^{(k)}_{i_1i_2\ldots i_m} > 0\}.$$ 
It follows that given any $i_2, \ldots, i_m$, there exists some $1 \le k \le K$ such that $p^{(k)}_{i_1i_2\ldots i_m} > 0$ for some $i_1 \in A$. According to  
\begin{equation}
\begin{split}
\Pr(\eta_{Ai_2\ldots i_m} \ge K+1) & =  \Pr\left(\bigcap_{j=m}^{m+K-1} \left\{X_j \notin A\right\} | X_{m-1}=i_2, \ldots, X_1=i_m\right)\\ \nonumber
 & = 1-\Pr\left(\bigcup_{j=m}^{m+K-1} \left\{X_j \in A\right\} | X_{m-1}=i_2, \ldots, X_1=i_m\right)\\
 & \le 1-\Pr(X_{m+k-1}=i_1 \in A | X_{m-1}=i_2, \ldots, X_1=i_m),
\end{split}
\end{equation}
where $1 \le k \le K$ is such that $p^{(k)}_{i_1i_2\ldots i_m}>0$, we see that 
\begin{equation}
\label{etaK}
\Pr(\eta_{Ai_2\ldots i_m} \ge K+1) \le 1-p^{(k)}_{i_1i_2\ldots i_m} < 1.
\end{equation}

On the other hand, 
$$\begin{array}{l}
\Pr(\eta_{Ai_2\ldots i_m} \ge K+1)\\
 =  \displaystyle \Pr\left(\bigcap_{j=m}^{m+K-1} \left\{X_j \notin A\right\} | X_{m-1}=i_2, \ldots, X_1=i_m\right)\\
 = \displaystyle \sum_{j_1 \in T}\sum_{j_2 \in T}\ldots \sum_{j_{K-1} \in T}\sum_{j_K \in T}p_{j_1i_2\ldots i_m}p_{j_2j_1i_2\ldots i_{m-1}}\cdots p_{j_{K-1}\ldots j_1i_2\ldots i_{m-K+2}}p_{j_K\ldots j_1i_2\ldots i_{m-K+1}},
\end{array}$$
where, and in what follows, if a multi-index is longer than $m$, then only the first $m$ indices actually count. This, together with (\ref{absoeqn}), i.e., $\tilde p^{(k)}_{i_1i_2\ldots i_m}=\delta_{i_1i_2}$ for any $i_2 \in A$, and (\ref{partpkeqn}), implies that
\begin{equation}
 \label{etaK1}
\begin{array}{l}
\Pr(\eta_{Ai_2\ldots i_m} \ge K+1)\\
 = \displaystyle \sum_{j_K \in T}\sum_{j_1 \in T}\ldots \sum_{j_{K-2} \in T}\sum_{j_{K-1} \in S} \tilde p_{j_Kj_{K-1}\ldots j_1i_2\ldots i_{m-K+1}}\tilde p_{j_{K-1}\ldots j_1i_2\ldots i_{m-K+2}}p_{j_{K-2}\ldots j_1i_2\ldots i_{m-K+3}}\cdots p_{j_1i_2\ldots i_m}\\
 = \displaystyle \sum_{j_K \in T}\sum_{j_1 \in T}\ldots \sum_{j_{K-3} \in T}\sum_{j_{K-2} \in S} \tilde p^{(2)}_{j_Kj_{K-2}\ldots j_1i_2\ldots i_{m-K+2}}\tilde p_{j_{K-2}\ldots j_1i_2\ldots i_{m-K+3}}p_{j_{K-3}\ldots j_1i_2\ldots i_{m-K+4}}\cdots p_{j_1i_2\ldots i_m}\\
 = \displaystyle \sum_{j_K \in T}\sum_{j_1 \in T}\ldots \sum_{j_{K-4} \in T}\sum_{j_{K-3} \in S} \tilde p^{(3)}_{j_Kj_{K-3}\ldots j_1i_2\ldots i_{m-K+3}}\tilde p_{j_{K-3}\ldots j_1i_2\ldots i_{m-K+4}}p_{j_{K-4}\ldots j_1i_2\ldots i_{m-K+5}}\cdots p_{j_1i_2\ldots i_m}\\
 \\
 = \ldots \\
 \\
 = \displaystyle \sum_{j_K \in T}\tilde p^{(K)}_{j_Ki_2\ldots i_m} = \sum_{i_1 \in T}\tilde p^{(K)}_{i_1i_2\ldots i_m}.
\end{array}
\end{equation}
The above expression, by (\ref{etaK}), gives $\displaystyle \sum_{i_1 \in T}\tilde p^{(K)}_{i_1i_2\ldots i_m} < 1$ for any $i_2, \ldots, i_m \in S$. Let 
$$\beta=\max_{i_2, \ldots, i_m \in S}\sum_{i_1 \in T}\tilde p^{(K)}_{i_1i_2\ldots i_m} < 1.$$

We now claim that there exists $0 \le \alpha < 1$ such that $\displaystyle \sum_{i_1 \in T}\tilde p^{(K+k)}_{i_1i_2\ldots i_m} < \alpha^{K+k}$ for any $i_2, \ldots, i_m \in S$ and $0 \le k \le K-1$. For $k=1$, 
$$\sum_{i_1 \in T}\tilde p^{(K+1)}_{i_1i_2\ldots i_m}=\sum_{i_1 \in T}\sum_{j \in S} \tilde p^{(K)}_{i_1ji_2\ldots i_{m-1}} \tilde p_{ji_2\ldots i_m} \le \beta\sum_{j \in S}\tilde p_{ji_2\ldots i_m}=\beta.$$
Similarly, for $k=2$, 
$$\sum_{i_1 \in T}\tilde p^{(K+2)}_{i_1i_2\ldots i_m}=\sum_{i_1 \in T}\sum_{j \in S} \tilde p^{(K+1)}_{i_1ji_2\ldots i_{m-1}} \tilde p_{ji_2\ldots i_m} \le \beta\sum_{j \in S}\tilde p_{ji_2\ldots i_m}=\beta.$$
Proceeding in this manner, we know that for any $i_2, \ldots, i_m \in S$ and $0 \le k \le K-1$, 
$$\sum_{i_1 \in T}\tilde p^{(K+k)}_{i_1i_2\ldots i_m} \le \beta.$$
Set 
$$\alpha=\beta^{\frac{1}{2K-1}}.$$
Observe that $0 \le \alpha < 1$. Then, for any $i_2, \ldots, i_m \in S$ and $0 \le k \le K-1$, 
$$\sum_{i_1 \in T}\tilde p^{(K+k)}_{i_1i_2\ldots i_m} \le \beta = \alpha^{2K-1} \le \alpha ^{K+k}.$$

Finally, we claim $\displaystyle \sum_{i_1 \in T}\tilde p^{(jK + k)}_{i_1i_2\ldots i_m} \le \alpha^{jK + k}$ for any $i_2, \ldots, i_m \in S$, $0 \le k \le K-1$, and $j=1, 2, \ldots$. This has been done for $j=1$. For $j=2$, using (\ref{absoeqn}) and (\ref{partpkeqn}) again and with the same technique as in the reduction of (\ref{etaK1}), we get 
\begin{equation}
\begin{split}
\sum_{i_1 \in T}\tilde p^{(2K+k)}_{i_1i_2\ldots i_m} & = \sum_{i_1 \in T}\sum_{j_1 \in T}\tilde p^{(2K+k-1)}_{i_1j_1i_2\ldots i_{m-1}}\tilde p_{j_1i_2\ldots i_m}\\ \nonumber
 & = \sum_{i_1 \in T}\sum_{j_1 \in T}\ldots \sum_{j_K \in T}\tilde p^{(K+k)}_{i_1j_K\ldots j_1i_2\ldots i_{m-K}}\tilde p_{j_K\ldots j_1i_2\ldots i_{m-K+1}}\cdots \tilde p_{j_1i_2\ldots i_m}\\
  & \le \alpha^{K+k}\sum_{j_K \in T}\sum_{j_1 \in T}\ldots \sum_{j_{K-1} \in T}\tilde p_{j_Kj_{K-1}\ldots j_1i_2\ldots i_{m-K+1}}\tilde p_{j_{K-1}\ldots j_1i_2\ldots i_{m-K+2}}\cdots \tilde p_{j_1i_2\ldots i_m}\\
  & = \alpha^{K+k}\sum_{j_K \in T}\sum_{j_1 \in T}\ldots \sum_{j_{K-2} \in T}\tilde p^{(2)}_{j_Kj_{K-2}\ldots j_1i_2\ldots i_{m-K+2}}\tilde p_{j_{K-2}\ldots j_1i_2\ldots i_{m-K+3}}\cdots \tilde p_{j_1i_2\ldots i_m}\\
  \\
  & = \ldots\\
  \\
  & = \alpha^{K+k}\sum_{j_K \in T}\tilde p^{(K)}_{j_Ki_2\ldots i_m} \le \alpha^{2K+k}.
\end{split}
\end{equation}
Similarly, for $j=3$, we have 
\begin{equation}
\begin{split}
\sum_{i_1 \in T}\tilde p^{(3K+k)}_{i_1i_2\ldots i_m} & = \sum_{i_1 \in T}\sum_{j_1 \in T}\tilde p^{(3K+k-1)}_{i_1j_1i_2\ldots i_{m-1}}\tilde p_{j_1i_2\ldots i_m}\\ \nonumber
 & = \sum_{i_1 \in T}\sum_{j_1 \in T}\ldots \sum_{j_K \in T}\tilde p^{(2K+k)}_{i_1j_K\ldots j_1i_2\ldots i_{m-K}}\tilde p_{j_K\ldots j_1i_2\ldots i_{m-K+1}}\cdots \tilde p_{j_1i_2\ldots i_m}\\
  & \le \alpha^{2K+k}\sum_{j_K \in T}\sum_{j_1 \in T}\ldots \sum_{j_{K-1} \in T}\tilde p_{j_Kj_{K-1}\ldots j_1i_2\ldots i_{m-K+1}}\tilde p_{j_{K-1}\ldots j_1i_2\ldots i_{m-K+2}}\cdots \tilde p_{j_1i_2\ldots i_m}\\
  & = \alpha^{2K+k}\sum_{j_K \in T}\sum_{j_1 \in T}\ldots \sum_{j_{K-2} \in T}\tilde p^{(2)}_{j_Kj_{K-2}\ldots j_1i_2\ldots i_{m-K+2}}\tilde p_{j_{K-2}\ldots j_1i_2\ldots i_{m-K+3}}\cdots \tilde p_{j_1i_2\ldots i_m}\\
  \\
  & = \ldots\\
  \\
  & = \alpha^{2K+k}\sum_{j_K \in T}\tilde p^{(K)}_{j_Ki_2\ldots i_m} \le \alpha^{3K+k}.
\end{split}
\end{equation}
Proceeding in this manner, the second claim can also be justified. The proof is now complete.
\ep

A bound similar to (\ref{tpk}) is developed in \cite{HX26b} using an alternative method and by resorting to (\ref{absolim}), see the proof of Corollary \ref{nonsing1} later for this bound. Theorem \ref{bound}, on the other hand, hinges solely on the ergodicity of the chain $X$, but not on (\ref{absolim}) at all. An immediate consequence of Theorem \ref{bound} is the following:

\begin{2}
\label{ergoabso}
Given the $(m-1)$th order ergodic chain $X$, let $\tilde X$ be the modified chain that is obtained from $X$ by making one or more of its states absorbing. Then, $\tilde X$ is an absorbing chain.
\end{2}
\bp
From (\ref{tpk}), we have $\ds \lim_{k \rightarrow \infty} \tilde p^{(k)}_{i_1i_2\ldots i_m}=0$ for any $i_1 \in T$ and $i_2, \ldots, i_m \in S$. The conclusion follows now from Definition \ref{abso}.
\ep

Incidentally, following the above theorem, (\ref{part}) is now the canonical form for the chain $\tilde X$.

In \cite{HX26b}, the fundamental tensor of an absorbing chain $\tilde X$ with transition tensor $\tilde {\cal P}$ in the form (\ref{part}) is defined by 
\be
\label{fund}
\sigma=\sum_{k=0}^\infty {\cal Q}^k={\cal I}+{\cal Q}+{\cal Q}^2+\ldots.
\ee
In addition, if it exists, $\sigma$ satisfies 
\be
\label{fundeqn}
\sigma={\cal I}+\sigma \boxtimes {\cal Q}.
\ee
Observe that for $m \ge 3$, the above equation may not be written as 
$$\sigma\boxtimes ({\cal I}-{\cal Q})={\cal I}$$
since $\cal I$ is only a left identity tensor.

Following (\ref{partpkeqn}) and Theorem \ref{bound}, we also obtain:
\begin{2}
\label{conv}
Given the $(m-1)$th order ergodic chain $X$, let $\tilde X$ be the modified chain obtained from $X$ by making $1 \le a < n$ states absorbing. We assume that the transition tensor of $\tilde X$ is given as in (\ref{part}). Then, the series in (\ref{fund}) converges and hence the fundamental tensor $\sigma$ is well defined.
\end{2}
\bp
The conclusion is obvious since ${\cal Q}=\tilde {\cal P}(T, T, :, \ldots, :)$ and, therefore, by Theorem \ref{bound},
$${\cal Q}^k \le \alpha^k {\cal E}$$
for any $k \ge K$, where $\cal E$ is the tensor of all ones of the same size as $\cal Q$.
\ep

For an $(m-1)$th order absorbing chain, the convergence of (\ref{fund}) follows from condition (\ref{absolim}) too \cite{HX26b}. The result above, however, is derived from the ergodicity of the chain $X$.

Also speaking of (\ref{fund}), it is well known \cite{KS} that for the case $m=2$, i.e., when $\tilde X$ is a first order absorbing chain, 
\be
\label{inv}
({\cal I}-{\cal Q})\boxtimes ({\cal I}+{\cal Q}+{\cal Q}^2+\ldots.) = ({\cal I}+{\cal Q}+{\cal Q}^2+\ldots.)\boxtimes ({\cal I}-{\cal Q}) = {\cal I}.
\ee
Thus, the fundamental tensor, i.e., now the fundamental matrix, is precisely the inverse of ${\cal I} -{\cal Q}$. Unfortunately, in general, (\ref{inv}) is no longer valid for a higher order absorbing chain, which is a noteworthy distinction between the first order and the higher order cases.

In addition, concerning (\ref{fundeqn}), we can state the following useful result.
\begin{2}
\label{nonsing}
Given the $(m-1)$th order ergodic chain $X$, let $\tilde X$ be the modified chain obtained from $X$ by making $1 \le a < n$ states absorbing. We assume that the transition tensor of $\tilde X$ has the form (\ref{part}). Then, as a linear system, (\ref{fundeqn}) is nonsingular, and its unique solution yields the fundamental tensor $\sigma$. 
\end{2}
\bp
First, (\ref{fundeqn}) is consistent since $\displaystyle \sigma=\sum_{k=0}^\infty {\cal Q}^k$ is a solution. Suppose that $\tilde \sigma=\sigma + \tau$ is another solution. Assume $\tau = [\tau_{i_1i_2\ldots i_m}] \ne \cal O$. By (\ref{fundeqn}), 
$$\tau=\tau \boxtimes {\cal Q},$$
and, therefore, it follows 
$$\tau=\underbrace{(\cdots ((\tau \boxtimes {\cal Q})\boxtimes {\cal Q})\boxtimes \cdots)\boxtimes {\cal Q}}_{k}$$
for any $k \ge K$. Next, on letting $\ds \rho=\max_{i_1, i_2 \in T; i_3, \ldots, i_m \in S}|\tau_{i_1i_2\ldots i_m}| > 0$ and making a similar argument as in (\ref{etaK1}), we obtain 
\begin{equation}
\begin{split}
|\tau_{i_1i_2\ldots i_m}| &= \big|\sum_{j_k \in T}\sum_{j_{k-1} \in T}\ldots \sum_{j_2 \in T}\sum_{j_1 \in T}\tau_{i_1j_1\ldots j_ki_2\ldots i_{m-k}}p_{j_1\ldots j_ki_2\ldots i_{m-k+1}}\cdots p_{j_{k-1}j_ki_2\ldots i_{m-1}}p_{j_ki_2\ldots i_m}\big|\\ \nonumber
 & \le  \rho \sum_{j_1 \in T}\sum_{j_k \in T}\ldots \sum_{j_3 \in T}\sum_{j_2 \in S} \tilde p_{j_1j_2\ldots j_ki_2\ldots i_{m-k+1}}\tilde p_{j_2\ldots j_ki_2\ldots i_{m-k+2}}\cdots p_{j_{k-1}j_ki_2\ldots i_{m-1}}p_{j_ki_2\ldots i_m}\\
 & =\rho \sum_{j_1 \in T}\sum_{j_k \in T}\ldots \sum_{j_4 \in T}\sum_{j_3 \in S} \tilde p^{(2)}_{j_1j_3\ldots j_ki_2\ldots i_{m-k+2}}\tilde p_{j_3\ldots j_ki_2\ldots i_{m-k+3}}\cdots p_{j_{k-1}j_ki_2\ldots i_{m-1}}p_{j_ki_2\ldots i_m}\\
 \\
 & = \ldots\\
 \\
 & = \rho \sum_{j_1 \in T} \tilde p^{(k)}_{j_1i_2\ldots i_m} \le \rho \alpha^k
\end{split}
\end{equation}
for any $i_1, i_2 \in T$ and $i_3, \ldots, i_m \in S$, implying 
$$\rho=\max_{i_1, i_2 \in T; i_3, \ldots, i_m \in S}|\tau_{i_1i_2\ldots i_m}| \le \rho \alpha^k.$$
This leads to a contradiction since $\alpha^k < 1$.
\ep

As mentioned earlier, a bound similar to (\ref{tpk}) is obtained for an $(m-1)$ order absorbing chain in \cite{HX26b} using (\ref{absolim}). Hence, the nonsingularity of (\ref{fundeqn}) and the uniqueness of the fundamental tensor as the solution of (\ref{fundeqn}) stated in Theorem \ref{nonsing} extend to such a chain. Specifically, we can state:

\begin{3}
\label{nonsing1}
Let $\tilde {\cal P}$ in (\ref{part}) be the transition tensor of an $(m-1)$th order absorbing chain, whose absorbing states are $A=\{t+1, t+2, \ldots, n\}$. Then, as a linear system, (\ref{fundeqn}) is nonsingular, and its unique solution yields the fundamental tensor. 
\end{3}
\bp
By the parallel bound in \cite{HX26b}, there exist $0<\alpha <1$ and $\beta>0$ such that for any $k \ge 1$, $i_1 \in T$, and $i_2, \ldots, i_m \in S$, 
$$\tilde p^{(k)}_{i_1i_2\ldots i_m} \le \beta \alpha^k.$$ 
The rest of the proof proceeds now in a similar fashion as that of Theorem \ref{nonsing}.
\ep

To determine the fundamental tensor $\sigma$ of an $(m-1)$th order absorbing chain, based on the preceding result, we can solve (\ref{fundeqn}) directly as a nonsingular linear system. This method has been implemented using MATLAB as {\tt fund} and added to the package described in \cite{Xu26c}, which is available at:

\centerline{\url{https://neumann.math.siu.edu/homc}}
\vspace{.2in}

To illustrate Theorems \ref{ergoabso}, \ref{conv}, and \ref{nonsing}, we give the following:

\begin{6}
\label{fundex1}
The transition tensor $\cal P$ of a second order chain on state space $S=\{1, 2, 3, 4, 5\}$ is given by 
$${\cal P}(:,:,1)=\left[\begin{array}{ccccc}
0 & 1/3 & 0 & 0 & 0\\
0 & 1/3 & 0 & 1/2 & 0\\
1 & 0 & 1 & 0 & 0\\
0 & 1/3 & 0 & 0 & 0\\
0 & 0 & 0 & 1/2 & 1
\end{array}\right], ~{\cal P}(:,:,2)=\left[\begin{array}{ccccc}
0 & 1/2 & 0 & 0 & 1/3\\
0 & 1/2 & 1/2 & 1 & 0\\
1 & 0 & 0 & 0 & 0\\
0 & 0 & 1/2 & 0 & 1/3\\
0 & 0 & 0 & 0 & 1/3
\end{array}\right],$$
$${\cal P}(:,:,3)=\left[\begin{array}{ccccc}
0 & 0 & 0 & 0 & 0\\
1/3 & 0 & 1/2 & 0 & 1\\
0 & 1 & 0 & 1/3 & 0\\
1/3 & 0 & 0 & 1/3 & 0\\
1/3 & 0 & 1/2 & 1/3 & 0
\end{array}\right], ~{\cal P}(:,:,4)=\left[\begin{array}{ccccc}
0 & 0 & 1/2 & 0 & 0\\
1 & 1 & 0 & 0 & 1/3\\
0 & 0 & 1/2 & 1 & 1/3\\
0 & 0 & 0 & 0 & 0\\
0 & 0 & 0 & 0 & 1/3
\end{array}\right],$$
$${\cal P}(:,:,5)=\left[\begin{array}{ccccc}
1/2 & 0 & 0 & 0 & 1/3\\
0 & 1/3 & 0 & 0 & 0\\
0 & 1/3 & 0 & 1/2 & 1/3\\
0 & 1/3 & 1 & 0 & 0\\
1/2 & 0 & 0 & 1/2 & 1/3
\end{array}\right].$$
Observe that ${\cal P}+{\cal P}^2+\ldots + {\cal P}^8 > 0$. This implies that for any $i_1, i_2, i_3 \in S$, there exists $1 \le k \le 8$, which may depend on $i_1, i_2, i_3$, such that $p^{(k)}_{i_1i_2i_3} > 0$. Consequently, the above chain is ergodic.

Next, we make, for example, states $2$ and $4$ absorbing. According to Theorem \ref{ergoabso}, this produces an absorbing chain. Denote by $\tilde {\cal P}$ its transition tensor. Use a simultaneous permutation over all the indices of $\tilde {\cal P}$ so as to bring state $2$ to the last position and $\tilde {\cal P}$ to the form in (\ref{part}), which can be done in MATLAB by \,{\tt x\,=[1 3 5 4 2]} \!and \,$\tilde {\tt P}$\,{\tt=}\,$\tilde {\tt P}${\tt (x,x,x)}. In particular, the ${\cal Q}$ block as in (\ref{part}) is given by 
$${\cal Q}(:,:,1)=\left[\begin{array}{ccc}
0 & 0 & 0\\
1 & 1 & 1\\
0 & 0 & 1
\end{array}\right], ~{\cal Q}(:,:,2)=\left[\begin{array}{ccc}
0 & 0 & 0\\
0 & 0 & 0\\
1/3 & 1/2 & 0
\end{array}\right], ~{\cal Q}(:,:,3)=\left[\begin{array}{ccc}
1/2 & 0 & 1/3\\
0 & 0 & 1/3\\
1/2 & 0 & 1/3
\end{array}\right],$$
$${\cal Q}(:,:,4)=\left[\begin{array}{ccc}
0 & 1/2 & 0\\
0 & 1/2 & 1/3\\
0 & 0 & 1/3
\end{array}\right], ~{\cal Q}(:,:,5)=\left[\begin{array}{ccc}
0 & 0 & 1/3\\
1 & 0 & 0\\
0 & 0 & 1/3
\end{array}\right].$$
Finally, we solve (\ref{fundeqn}) using {\rm \textsigma}{\tt\,=\,fund(Q)} to obtain the fundamental tensor:
$$\sigma(:,:,1)=\left[\begin{array}{ccc}
1 & 0 & 1\\
2 & 2 & 4/3\\
1/2 & 1/2 & 7/2
\end{array}\right], ~\sigma(:,:,2)=\left[\begin{array}{ccc}
4/3 & 0 & 0\\
4/9 & 1 & 0\\
7/6 & 1/2 & 1
\end{array}\right], ~\sigma(:,:,3)=\left[\begin{array}{ccc}
2 & 0 & 1\\
10/6 & 1 & 4/3\\
2 & 0 & 5/2
\end{array}\right],$$
$$\sigma(:,:,4)=\left[\begin{array}{ccc}
1 & 2/3 & 1/3\\
0 & 31/18 & 7/9\\
0 & 5/6 & 11/6
\end{array}\right], ~\sigma(:,:,5)=\left[\begin{array}{ccc}
1 & 0 & 1\\
2 & 1 & 1\\
1/2 & 0 & 5/2
\end{array}\right].$$
\end{6}

As pointed out in the introductory section, results on higher order chains do not stem from the corresponding ones on first order chains in general. On the other hand, results on higher order chains do include the corresponding ones on first order chains as their special cases. Here is one such example to illustrate Corollary \ref{nonsing1}.

\begin{6}
\label{fundex2}
Take a first order chain on state space $S=\{1, 2, 3, 4, 5, 6\}$ with transition matrix 
$$\tilde {\cal P}=\left[\begin{array}{cccccc}
0 & 1/2 & 0 & 0 & 0 & 0\\
1 & 0 & 1/3 & 0 & 0 & 0\\
0 & 1/2 & 1/3 & 0 & 0 & 0\\
0 & 0 & 0 & 1 & 0 & 0\\
0 & 0 & 1/3 & 0 & 1 & 0\\
0 & 0 & 0 & 0 & 0 & 1
\end{array}\right].$$
The matrix $\tilde {\cal P}$ is already in the form (\ref{part}). It is easy to see $\ds \lim_{k \rightarrow \infty}{\cal Q}^k = {\cal O}$, the zero tensor of the same size as $\cal Q$, thus the above chain is absorbing. Note   
$${\cal Q}=\left[\begin{array}{ccc}
0 & 1/2 & 0\\
1 & 0 & 1/3\\
0 & 1/2 & 1/3
\end{array}\right].$$
With {\rm \textsigma}{\tt\,=\,fund(Q)}, the fundamental matrix $\sigma$ is found to be 
$$\sigma=\left[\begin{array}{ccc}
3 & 2 & 1\\
4 & 4 & 2\\
3 & 3 & 3
\end{array}\right],$$
which coincides with $({\cal I}-{\cal Q})^{-1}$.
\end{6}

Our next question concerns the relationship between the mean first passage tensor $\mu$ of the $(m-1)$th order ergodic chain $X$ and the fundamental tensor $\sigma$ of the modified absorbing chain $\tilde X$. In what follows, we shall focus on the case where only one of the states is made absorbing. Without loss of generality, assume that it is state $n$, i.e., $T=\{1, 2, \ldots, n-1\}$ and $A=\{n\}$ in what follows.

Observe that the transition tensor of the ergodic chain $X$ is partitioned as in (\ref{part0}). From (\ref{mfpteqn}), the mean first passage times on $X$ to state $n$ are the unique solutions to 
\be
\label{mfpteqn1}
\mu_{ni_2\ldots i_m}=1+\sum_{j \in T} \mu_{nji_2\ldots i_{m-1}}p_{ji_2\ldots i_m}, ~i_2 \in T,
\ee
and 
\be
\label{mfpteqn2}
\mu_{nni_3\ldots i_m}=1+\sum_{j \in T} \mu_{njni_3\ldots i_{m-1}}p_{jni_3\ldots i_m},
\ee
for any $i_3, \ldots, i_m \in S$. Clearly, all these times form a horizontal slice $\mu(n,:,\ldots,:)$ of $\mu$. On letting 
$$\mu(n,:,\ldots,:)=[\nu ~~\omega],$$
where $\nu$ and $\omega$ are $m$th order tensors, $\nu=\mu(n,T,:,\ldots,:) \in \mathbb R^{1 \times (n-1) \times n \times \cdots \times n}$ and $\omega=\mu(n,n,:,\ldots,:) \in \mathbb R^{1 \times 1 \times n \times \cdots \times n}$, (\ref{mfpteqn1}) and (\ref{mfpteqn2}) can be written, in the tensor form, more conveniently as 
\be
\label{nu}
\nu={\cal E}+\nu \boxtimes {\cal Q},
\ee
where $\cal E$ is the tensor of all ones and of the same size as $\nu$, and
\be
\label{ome}
\omega={\cal E}+\rho \boxtimes {\cal S},
\ee
where $\cal E$ is the tensor of all ones and of the same size as $\omega$, and where $\rho$ is the tensor of the same size as $\nu$ such that 
\be
\label{rho}
\rho(:,:,i_3,:,\ldots,:)=\nu(:,:,n,:,\ldots,:)
\ee
for any $i_3 \in S$. According to (\ref{nu}) and (\ref{ome}), to determine $\mu(n,:,\ldots,:)$, we can first solve (\ref{nu}) for $\nu$, and then use (\ref{ome}) to find $\omega$. 

Before proceeding, let us state the following last technical lemma.

\begin{1}
Let ${\cal A}=[a_{i_1i_2\ldots i_m}] \in \mathbb R^{u \times w \times n \times \cdots \times n}$, where $w \le n$, be an $m$th order tensor and ${\cal E} \in \mathbb R^{1 \times u \times n \times \cdots \times n}$ be an $m$th order tensor of all ones. Then, ${\cal E} \boxtimes {\cal A} \in \mathbb R^{1 \times w \times n \times \cdots \times n}$ is an $m$th order tensor such that 
$${\cal E} \boxtimes {\cal A}=\left[\sum_{j=1}^u a_{ji_2\ldots i_m}\right],$$
i.e., it effects a summation operation over the first index of $\cal A$. In particular, for an $m$th order identity tensor ${\cal I} \in \mathbb R^{u \times u \times n \times \cdots \times n}$, where $u \le n$, 
\be
\label{eeye}
{\cal E} \boxtimes {\cal I}={\cal E}.
\ee
Moreover, if ${\cal B} \in \mathbb R^{w \times v \times n \times \cdots \times n}$, where $v \le n$, is an $m$th order tensor too, then 
\be
\label{asso}
({\cal E} \boxtimes {\cal A}) \boxtimes {\cal B}={\cal E} \boxtimes ({\cal A} \boxtimes {\cal B}),
\ee
i.e., there is associativity for such a special case of $\boxtimes$ multiplication.
\end{1}
\bp
The results regarding ${\cal E} \boxtimes {\cal A}$ and ${\cal E} \times {\cal I}$ are straightforward. To show (\ref{asso}), we calculate $({\cal E} \boxtimes {\cal A})\boxtimes {\cal B}$ to get 
\be
\begin{split}
({\cal E} \boxtimes {\cal A})\boxtimes {\cal B} & =\left[\sum_{j=1}^w \sum_{i_1=1}^u a_{i_1ji_2\ldots i_{m-1}}b_{ji_2\ldots i_m}\right]=\left[\sum_{i_1=1}^u \sum_{j=1}^w a_{i_1ji_2\ldots i_{m-1}}b_{ji_2\ldots i_m}\right]\\ \nonumber
 & = {\cal E} \boxtimes ({\cal A} \boxtimes {\cal B}).
\end{split}
\ee
\ep

The next result establishes the relationship between the component $\nu$ of $\mu(n,:,\ldots,:)$ and the fundamental tensor $\sigma$. It also gives a series representation for $\nu$.

\begin{2}
\label{mfptfund}
For the $(m-1)$th order ergodic chain $X$ whose transition tensor ${\cal P}$ is partitioned as (\ref{part0}), let $\tilde X$ be the modified absorbing chain whose transition tensor has the form (\ref{part}) with $A=\{n\}$.  Then, $\nu$ in (\ref{nu}) can be represented by 
\be
\label{nuser}
\nu={\cal E} \boxtimes \sigma = \sum_{k=0}^\infty  {\cal E} \boxtimes {\cal Q}^k,.
\ee
where $\sigma$ is the fundamental tensor of $\tilde X$ and ${\cal E} \in \mathbb R^{1 \times (n-1) \times n \times \cdots \times n}$ is an $m$th order tensor of all ones. In addition, as a linear system, (\ref{nu}) is nonsingular, and its unique solution yields $\nu$ in (\ref{nuser}). Once $\nu$ is available, $\omega$ can be determined from (\ref{ome}).
\end{2}
\bp
According to Theorem \ref{conv}, $\ds \sigma=\sum_{k=0}^\infty {\cal Q}^k$ is well defined.

First, we show that $\nu$ in (\ref{nuser}) satisfies (\ref{nu}). This is easy to check since using (\ref{eeye}) and (\ref{asso}), we see 
\be
\begin{split}
{\cal E}+\nu \boxtimes {\cal Q} &= {\cal E} \boxtimes {\cal I}+({\cal E} \boxtimes ({\cal I} + {\cal Q} + {\cal Q}^2 + \ldots ))\boxtimes {\cal Q}\\ \nonumber
 & = {\cal E} \boxtimes ({\cal I} + {\cal Q} + {\cal Q}^2 + \ldots)=\nu.
 \end{split}
\ee
 Next, the nonsingularity of (\ref{nu}) can be verified in a similar manner as the proof of Theorem \ref{nonsing}, which we shall omit for that reason.
\ep

To illustrate Theorem \ref{mfptfund}, let us consider:
\begin{6}
\label{mfptex1}
Take a second order chain on $S=\{1, 2, 3, 4\}$ whose transition tensor is given by 
$${\cal P}(:,:,1)=\left[\begin{array}{cccc}
1/2 & 0 & 0 & 0\\
1/2 & 0 & 1 & 0\\
0 & 1 & 0 & 1\\
0 & 0 & 0 & 0
\end{array}\right], ~{\cal P}(:,:,2)=\left[\begin{array}{cccc}
0 & 0 & 1/2 & 1\\
0 & 1/2 & 0 & 0\\
0 & 1/2 & 0 & 0\\
1 & 0 & 1/2 & 0
\end{array}\right],$$
$${\cal P}(:,:,3)=\left[\begin{array}{cccc}
0 & 1 & 0 & 1\\
1 & 0 & 1/2 & 0\\
0 & 0 & 1/2 & 0\\
0 & 0 & 0 & 0
\end{array}\right], ~{\cal P}(:,:,4)=\left[\begin{array}{cccc}
0 & 0 & 0 & 0\\
1 & 1 & 1 & 0\\
0 & 0 & 0 & 1/2\\
0 & 0 & 0 & 1/2
\end{array}\right].$$
This chain is ergodic since ${\cal P}+{\cal P}^2+{\cal P}^3+{\cal P}^4 > {\cal O}$.

Now, we use Theorem \ref{mfptfund} to calculate the mean first passage times of the above chain to, for example, state $2$. Like in Example \ref{fundex1}, $\cal P$ can be converted to the form (\ref{part}) via a simultaneous permutation, i.e., using {\tt x\,=[1 3 4 2]}, over all its indices. In particular, the $\cal Q$ and $\cal S$ blocks are given by 
$${\cal Q}(:,:,1)=\left[\begin{array}{ccc}
1/2 & 0 & 0\\
0 & 0 & 1\\
0 & 0 & 0
\end{array}\right], ~{\cal Q}(:,:,2)=\left[\begin{array}{ccc}
0 & 0 & 1\\
0 & 1/2 & 0\\
0 & 0 & 0
\end{array}\right],$$
$${\cal Q}(:,:,3)=\left[\begin{array}{ccc}
0 & 0 & 0\\
0 & 0 & 1/2\\
0 & 0 & 1/2
\end{array}\right], ~{\cal Q}(:,:,4)=\left[\begin{array}{ccc}
0 & 1/2 & 1\\
0 & 0 & 0\\
1 & 1/2 & 0
\end{array}\right],$$
and 
$${\cal S}(:,:,1)=[0 \ \ 1 \ \ 0]^T, ~{\cal S}(:,:,2)=[1 \ \ 0 \ \ 0]^T, ~{\cal S}(:,:,3)=[0 \ \ 0 \ \ 0]^T, ~{\cal S}(:,:,4)=[0 \ \ 1/2 \ \ 0]^T.$$
According to Theorem \ref{ergoabso}, the modified chain $\tilde X$ having transition tensor in the form of (\ref{part}) is absorbing. Incidentally, there is no need to construct the transition tensor of $\tilde X$. All we need here are the $\cal Q$ and $\cal S$ blocks.

Next, using {\rm \textsigma}{\tt\,=\,fund(Q)}, the fundamental tensor $\sigma$ turns out to be 
$$\sigma(:,:,1)=\left[\begin{array}{ccc}
2 & 0 & 0\\
0 & 1 & 1\\
0 & 0 & 1
\end{array}\right], ~\sigma(:,:,2)=\left[\begin{array}{ccc}
1 & 0 & 1\\
0 & 2 & 0\\
0 & 0 & 1
\end{array}\right],$$
$$\sigma(:,:,3)=\left[\begin{array}{ccc}
1 & 0 & 0\\
0 & 1 & 1\\
0 & 0 & 2
\end{array}\right], ~\sigma(:,:,4)=\left[\begin{array}{ccc}
1 & 1 & 1\\
1 & 1 & 0\\
1 & 1/2 & 1
\end{array}\right].$$
By (\ref{nuser}), we use {\rm \textnu}{\tt \,=\,bprod({\rm \textEpsilon},{\rm \textsigma})}, where {\rm \textEpsilon}{\tt\,=\,ones(1,3,4)}, to get
$$\nu(:,:,1)=[2 \ \ 1 \ \ 2], ~\nu(:,:,2)=[1 \ \ 2 \ \ 2], ~\nu(:,:,3)=[1 \ \ 1 \ \ 3], ~\nu(:,:,4)=[3 \ \ 5/2 \ \ 2].$$
It follows 
$$\rho(:,:,1)=\rho(:,:,2)=\rho(:,:,3)=\rho(:,:,4)=[3 \ \ 5/2 \ \ 2].$$
Note that when $m > 3$, $\rho(:,:,i_3,:,\ldots,:)$ in (\ref{rho}) may not equal each other. Moving on, from (\ref{ome}), we see {\rm \textomega}{\tt\,=\,}{\rm \textEpsilon}{\tt\,+\,bprod({\rm \textrho},S)}, where {\rm \textEpsilon}{\tt\,=\,ones(1,1,4)}. Hence,
$$\omega(:,:,1)=7/2, ~\omega(:,:,2)=4, ~\omega(:,:,3)=1, ~\omega(:,:,4)=9/2.$$
Applying now the reverse simultaneous permutation over the last two indices of  $\xi=[\nu ~~\omega]$, i.e., {\rm \textxi}{\tt \,=\,}{\rm \textxi}{\tt (:,y,y)} with {\tt y\,=[1 4 2 3]}, we obtain the horizontal slice $\mu(2,:,:)$ as:
$$\mu(2,:,1)=[2 \ \ 7/2 \ \ 1 \ \ 2], ~\mu(2,:,2)=[3 \ \ 9/2 \ \ 2.5 \ \ 2],$$
$$\mu(2,:,3)=[1 \ \ 4 \ \ 2 \ \ 2], ~\mu(2,:,4)=[1 \ \ 1 \ \ 1 \ \ 3].$$
\end{6}

In the same spirit as Example \ref{fundex2}, we also provide below an example for Theorem \ref{mfptfund} in the first order case.

\begin{6}
\label{mfptex2}
Let $X$ be a first order chain with transition matrix 
$${\cal P}=\left[\begin{array}{cccccc}
0 & 1/2 & 0 & 0 & 0 & 0\\
1 & 0 & 1/3 & 1/2 & 0 & 0\\
0 & 1/2 & 1/3 & 0 & 0 & 1/3\\
0 & 0 & 0 & 0 & 1 & 0\\
0 & 0 & 1/3 & 0 & 0 & 1/3\\
0 & 0 & 0 & 1/2 & 0 & 1/3
\end{array}\right]$$
and state space $S=\{1, 2, 3, 4, 5, 6\}$. Then $X$ is ergodic since ${\cal P}+{\cal P}^2+{\cal P}^3+{\cal P}^4+{\cal P}^5 > {\cal O}$.

To find, for example, the mean first passage times of $X$ to state $6$, we set ${\cal Q}={\cal P}(1\!:\!5,1\!:\!5)$. Using {\rm \textsigma}{\tt\,=\,fund(Q)}, we get the fundamental matrix $\sigma$ as:
$$\sigma=\left[\begin{array}{ccccc}
5 & 4 & 3 & 2 & 2\\
8 & 8 & 6 & 4 & 4\\
6 & 6 & 6 & 3 & 3\\
2 & 2 & 2 & 2 & 2\\
2 & 2 & 2 & 1 & 2
\end{array}\right].$$
Hence, by (\ref{nuser}) and (\ref{ome}), we calculate {\rm \textnu}{\tt \,=\,bprod({\rm \textEpsilon},{\rm \textsigma})}, where {\rm \textEpsilon}{\tt\,=\,ones(1,5)}, to obtain  
$$\nu=[23 ~~22 ~~19 ~~12 ~~13],$$ 
and {\rm \textomega}{\tt\,=\,1\,+\,bprod({\rm \textnu},S)} to find  
$$\omega=35/3.$$ 
Notice that for the latter, $\rho=\nu$ as in any first order case. Subsequently, the desired mean first passage times are given by  
$$\mu(6,:)=[23 ~~22 ~~19 ~~12 ~~13 ~~35/3].$$
\end{6}

Incidentally, the results in Examples \ref{mfptex1} and \ref{mfptex2} are consistent with the results computed with {\tt mfptd(P)} and {\tt mfpti(P)} in \cite{Xu26c}.

\section{Concluding Remarks}
\label{concl}
\setcounter{equation}{0}

In this paper, we further investigate the mean first passage times of a higher order ergodic chain in a setting when some states of the chain are modified to be absorbing. Our main contributions include: (1) Such a modified higher order chain is absorbing. (2) For a higher order absorbing chain, the tensor equation its fundamental tensor satisfies is a nonsingular linear system. (3) The mean first passage time tensor of a higher order ergodic chain is closely related to the fundamental tensor of a certain higher order absorbing chain. This also renders a new way of computing selected mean first passage times.

The MATLAB function {\tt fund} for computing the fundamental tensor of a higher order absorbing chain incorporates the direct method for solving the nonsingular linear system mentioned above, i.e.,
$$\sigma={\cal I}+\sigma \boxtimes {\cal Q}.$$
This system may also be solved in an iterative manner. Starting with $\sigma^{(0)}={\cal I}$, for example, we have 
$$\sigma^{(k+1)}={\cal I}+\sigma^{(k)} \boxtimes {\cal Q}, ~k=0, 1, 2, \ldots.$$
This approach, however, appears brute force since it actually computes the successive partial sums of the series $\ds \sigma=\sum_{k=0}^\infty {\cal Q}^k$. This is the reason why only the direct method is adopted. There may be more efficient iterative schemes for computing the fundamental tensor, which will be one area for our future work.

The regularity of higher order chains has been introduced in \cite{HX26a}. Specifically, a higher order chain with transition tensor $\cal P$ is regular if there exists $k \ge 1$ such that ${\cal P}^k > {\cal O}$. For the first order case, an alternative fundamental tensor, i.e., matrix in this case, can also be defined for a regular chain, which plays a significant role in the study of such a chain \cite{Ios, KS}. This certainly raises the question of whether similar results may be developed for a higher order chain, leading to another area for our future work.

\end{document}